\documentclass[9pt,
]{article}

  \usepackage[utf8]{inputenc}
  \usepackage{amsmath}
\usepackage{enumerate}
  \usepackage{amssymb}
  \usepackage{amsthm}
  \usepackage{bm}
  \usepackage{graphicx}
  \usepackage{color}
\usepackage{accents}

\newtheorem{thank}{\ \ \ Acknowledgment}

\def\scalar(#1,#2){(#1\mid#2)}

\newcommand{\ca}{\mathcal{A}}

\newcommand{\Pro}{{\mathbb{P}}}
\newcommand{\T}{{\mathbb{T}}}

\newcommand{\N}{{\mathbb{N}}}

\newcommand{\tend}[3][]{\xrightarrow[#2\to#3]{#1}}

\newcommand{\ds}{\displaystyle}

\title{A simple proof of Bourgain's theorem on the singularity of the spectrum of Ornstein's maps
 {\footnote{The reader need not be familiar with Ornstein's construction neither with the spectral theory of
  dynamical systems.}}}
\author{E. H. el Abdalaoui\footnote{Normandy University, University of Rouen
Department of Mathematics, LMRS  UMR 60 85 CNRS,
Avenue de l'Universit\'e, BP.12,
76801 Saint Etienne du Rouvray - France .}}
\date{October 2, 2014}
\begin{document}
\maketitle
{\renewcommand\abstractname{{\textbf{Abstract}}}
\abstractname. We give a simple proof of Bourgain's theorem on the singularity of Ornstein's maps.
\paragraph{Setting and proof.}
Let $(m_j)$, $(t_j)$ be a sequence of positive integers, and let $(\Omega,\ca,\Pro)$ be the probability space associated to Ornstein's construction, that is,
$$\Omega=\prod_{j=1}^{+\infty}\big\{-t_j,\cdots,t_j\big\}^{p_j-1},~~
\Pro=\bigotimes_{j=1}^{+\infty}\bigotimes_{k=1}^{p_j-1}{\mathcal{U}}_k,$$
where ${\mathcal{U}}_k$ is the uniform measure on $\big\{-t_j,\cdots,t_j\big\}^{p_j-1}.$\\

We want to prove the following theorem due to Bourgain \cite{BourgainD}.\\

{\it For almost all $\omega \in \Omega$, the spectral type $\mu_{\omega}$ of the rank one map $T_{\omega}$ is singular.}\\

 Recall that $\mu_{\omega}$ is the weak-star limit of the following sequence of probability measures
$$  \Big(\prod_{j=1}^{N}|P_j(\omega,z)|^2 d\lambda\Big)_{N \geq 1},$$
where $\lambda$ is the Lebesgue measure and for each $j \in \N^*$,
$$P_j(z)=\frac1{\sqrt{m_j}}\sum_{k=0}^{m_j-1}z^{n_{j,k}(\omega)},$$
$$n_{j,0}(\omega)=0~~~~{\rm{and~~~for~~}} k \geq 1,
n_{j,k}(\omega)=k(h_j+t_j)+x_{j,k}(\omega)=k(h_j+t_j)+\omega_{j,k}.$$
We further assume that the sequence $(m_j)$ is unbounded. Therefore,
 by Theorem 5.2 in \cite{Abd-Nad}.
 combined with the uniform integrability of the sequence $\prod_{j=1}^{N}|P_j(\omega,z)|$, we have
$$\int_{\Omega}\int \prod_{j=1}^{N}|P_j(\omega,z)| dz d\Pro \tend{N}{+\infty}
\int_{\Omega}\int \sqrt{\frac{d\mu_{\omega}}{d\lambda}}dzd\Pro,$$
We further have
$$\lim_{j \longrightarrow +\infty}\int |P_j(\omega,z)| d\Pro dz=\lim_{j \longrightarrow +\infty}\int |\widetilde{P}_j(\omega,z)| d\Pro dz =\frac12 \sqrt{\pi}, {~~\textrm{with}~~} \widetilde{P}_j(\omega,z)=P_j(\omega,z)-\int P_j(\omega,z)d\Pro. $$
 The last equality follows from the classical  Lindeberg's central limit theorem (CLT) \footnote{It is an easy exercise to check that the Lindeberg condition holds under $dz \otimes d\Pro$.} combined with Lebesgue dominated convergence theorem and the uniform integrability of the sequence $(|P_j(\omega,z)|)_{j \geq 0}$ under $dz \otimes \Pro$.
We thus get
$$\int \int_{\Omega} \prod_{j=1}^{N}|P_j(\omega,z)| d\Pro dz=\int \prod_{j=1}^{N}\int_{\Omega} |P_j(\omega,z)| d\Pro dz    \tend{N}{+\infty} 0.$$
Whence
$$\int_{\Omega}\int \sqrt{\frac{d\mu_{\omega}}{d\lambda}}dzd\Pro=0,$$
and the proof is complete.\qed
\section*{$\bigstar$More details.}Notice that
\begin{eqnarray*}
\int_{\Omega}\int_{\T}\Big||P_j(\omega,z)|-|\widetilde{P}_j(\omega,z)|\Big|dz d\Pro
&\leq& \int_{\T} \Big|\int P_m(\omega,z)d\Pro\Big|dz =
\int_{\T} \Big|\frac1{\sqrt{m}}\sum_{p=0}^{m-1}z^{p(h_m+t_m)}\Big| \Big|\frac1{t_m+1}\sum_{s=0}^{t_m}z^s\Big|dz\\
&\leq& \Big\|\frac1{\sqrt{m}}\sum_{p=0}^{m-1}z^{p(h_m+t_m)}\Big\|_2 \Big\|\frac1{t_m+1}\sum_{s=0}^{t_m}z^s\Big\|_2.
\end{eqnarray*}
The last inequality is due to the Cauchy-Schwarz inequality. This gives
\begin{eqnarray*}
\int_{\Omega}\int_{\T}\Big||P_j(\omega,z)|-|\widetilde{P}_j(\omega,z)|\Big|dz d\Pro
\leq \frac1{\sqrt{t_m+1}} \tend{m}{+\infty}0.
\end{eqnarray*}
Since
$$\Big\|\frac1{\sqrt{m}}\sum_{p=0}^{m-1}z^{p(h_m+t_m)}\Big\|_2=
\Big\|\frac1{\sqrt{m}}\sum_{p=0}^{m-1}z^{p}\Big\|_2=1, \textrm{~~and~~}
\Big\|\frac1{t_m+1}\sum_{s=0}^{t_m}z^s\Big\|_2=\frac1{\sqrt{t_m+1}}.$$
\section*{$\blacktriangleright$ On the proof of Theorem 5.2. in  ``Calculus of Generalized Riesz Products" pages 158-162.}
The proof of Theorem 5.2. is self-contained and goes as Follows.
\begin{itemize}
  \item Using Cauchy-Schwarz inequality, we establish that $\ds \sqrt{\frac{d\mu}{d\lambda}}$ is a weak limit of
  the sequence of $L^2$-functions $\ds \prod_{j=1}^{n}|P_j(z)|$.
  \item  We take advantage  of the following formula
  $$\mu=R_n^2(z) d\mu_n,$$
  where
  $$R_n=\prod_{j=1}^{n}|P_j(z)|~~~~~~~~~~~~~{\textrm{and}}~~~~~~~~~~ d\mu_n=\prod_{j=n+1}^{+\infty}|P_j(z)|^2,
  $$ and we prove that any weak limit $\phi$ of the sequence  $\ds \sqrt{\frac{d\mu_n}{d\lambda}}$ satisfy  $0 \leq \phi \leq 1$. Finally, we deduce that the limit of the sequence
  $ \Big\|R_n-\sqrt{\frac{d\mu}{d\lambda}}\Big\|_1$ is zero.
\end{itemize}
\section*{$\blacktriangleright$ On the CLT argument.}
It is easy to see that Lindeberg's condition holds for the sequence of the random variables
$$X_m(\omega,z)=\frac1{\sqrt{p_m}}\sum_{k=0}^{p_m-1}\Big(z^{n_{j,k}(\omega)}-\int_{\Omega}z^{n_{j,k}(\omega)}d\Pro\Big).$$
We thus get that $(X_m(\omega,z))$ converge in distribution to the complex normal distribution under $dz \otimes d\Pro$.
\begin{thank}
The author wishes to express his thanks to Fran\c cois Parreau, Jean-Paul Thouvenot, Mahendra Nadkarni and Mariusz Lema\'{n}czyk for discussions on this note.\\
It is also a pleasure for him to express his thanks for the warm hospitality
to the organizers of the conference
''Number Theory and Dynamics, 06 July - 10 July 2015, institut Mittag-Leffler", and to
the institut of Mittag-Leffer where this note was revised and completed.
\end{thank}

\end{document}